\documentclass[12pt,reqno]{amsart} 

\usepackage{mathrsfs,amsfonts,amsthm,amsmath,amssymb,thmtools} 
\usepackage[dvipsnames]{xcolor} 

\usepackage{enumitem}                
\usepackage{lineno}                 

\usepackage[numbers,sort&compress]{natbib}  

\usepackage{graphics}
\usepackage{graphicx}   
\graphicspath{{images/}} 
\usepackage{subfig}     
\usepackage{rotating}

\usepackage[hyperindex, pdfencoding = auto, psdextra, bookmarksdepth = 4]{hyperref}  
\hypersetup{
    bookmarksopen      = true,
    bookmarksopenlevel = 1,
    bookmarksnumbered  = true,
    pdftitle    = {},
    pdfauthor   = {},
    linktoc     = page,
    anchorcolor = green,
    colorlinks  = true,        
    linkcolor   = black,
    citecolor   = blue,
    filecolor   = blue,
    urlcolor    = magenta,
}

\usepackage[capitalise, noabbrev, nameinlink]{cleveref}    

\theoremstyle{plain}

\newtheorem{theorem}{Theorem}[section]
\newtheorem{lemma}[theorem]{Lemma}

\newtheorem{proposition}[theorem]{Proposition}
\newtheorem{corollary}[theorem]{Corollary}

\theoremstyle{definition}   

\makeatletter
    \def\subsection{\@startsection{subsection}{2}%
    \z@{.5\linespacing\@plus.7\linespacing}{.3\linespacing}%
    {\normalfont\bfseries}}
\makeatother 

\makeatletter                                         
    \newcommand{\LeftEqNo}{\let\veqno\@@leqno}        
\makeatother

\numberwithin{equation}{section}                      

\begin{document}

\

\vspace{-2cm}

\title[Embeddings of ultraproducts into group von Neumann algebras]
{Embeddings of ultraproducts into group von Neumann algebras}

\author{Minghui Ma}
\address{Minghui Ma, School of Mathematical Sciences, Dalian University of Technology, Dalian, 116024, China}
\email{minghuima@dlut.edu.cn}

\keywords{Ultraproduct, finite factor, Jones index, group von Neumann algebra, Fourier height}

\makeatletter
\@namedef{subjclassname@2020}{\textup{2020} Mathematics Subject Classification}
\makeatother
\subjclass[2020]{Primary 46L10; Secondary 46L37, 46M07}

\begin{abstract}
Let $(\mathcal M_n,\tau_n)_{n\geqslant1}$ be a sequence of finite factors and let $\omega$ be a free ultrafilter on $\mathbb N$ such that $\dim_{\mathbb C}\mathcal M_n\to_\omega\infty$.
We prove that every subfactor of finite Jones index in the ultraproduct $\prod_\omega\mathcal M_n$ admits no embedding into the group von Neumann algebra of any discrete group.
In particular, neither the ultrapower $\mathcal R^\omega$ of the hyperfinite type $\mathrm{II}_1$ factor nor the matrix ultraproduct $\prod_\omega M_n(\mathbb C)$ embeds into the group von Neumann algebra of any discrete group.
\end{abstract}

\maketitle

\section{Introduction}

Group von Neumann algebras provide a fundamental class of examples of type $\mathrm{II}_1$ factors, dating back to the work of Murray and von Neumann \cite{MvN-43}.
These include the hyperfinite type $\mathrm{II}_1$ factor $\mathcal R$ and the free group factor $\mathcal L(\mathbb F_2)$.
Sakai \cite[Problem~4.4.30]{Sakai-71} recorded the question whether every type $\mathrm{II}_1$ factor is a group factor.
Connes \cite[Corollary~7]{Connes-75} answered this question in the negative by constructing a type $\mathrm{II}_1$ factor that is not isomorphic to its opposite algebra.
Jones \cite{Jones-80} subsequently constructed a nongroup factor isomorphic to its opposite but admitting no involutory antiautomorphism.
Later, Ioana \cite[Corollary~10.1 and Remark~10.3(3)]{Ioana-11} obtained nongroup factors admitting involutory antiautomorphisms, while Boutonnet \cite[Theorem~4.1]{Boutonnet-13} constructed factors not stably isomorphic to any group factor.

Which factors can occur as subfactors of group factors?
Ozawa \cite[Corollary~3]{Ozawa-04} showed that no single separable type $\mathrm{II}_1$ factor contains a copy of every separable type $\mathrm{II}_1$ factor.
Chifan, Drimbe, and Ioana \cite[Theorem~A and Remark~1.1]{CDI-23} proved that every separable tracial von Neumann algebra embeds into a type $\mathrm{II}_1$ factor with property~$(T)$ and asked whether every separable type $\mathrm{II}_1$ factor embeds into $\mathcal L(G)$ for some countable ICC group $G$.
In this paper, we consider this embedding problem for ultraproducts of finite factors and their subfactors of finite Jones index, allowing $G$ to be any discrete group.
Throughout, an \emph{embedding} means an injective unital $*$-homomorphism.

For a tracial von Neumann algebra $(\mathcal M,\tau)$ and $x,y\in\mathcal M$, we use the inner product $\langle x,y\rangle=\tau(xy^*)$ and the associated $2$-norm $\|x\|_{2,\tau}=\tau(x^*x)^{1/2}$.
The \emph{ultraproduct} of a sequence $(\mathcal M_n,\tau_n)_{n\geqslant 1}$ of finite factors along a free ultrafilter $\omega$ on $\mathbb N$ is defined by
\begin{equation*}
	\prod_\omega\mathcal M_n=\prod_{n=1}^\infty\mathcal M_n\bigg/\mathcal I_\omega,\qquad
	\mathcal I_\omega=\left\{(x_n)\in\prod_{n=1}^\infty\mathcal M_n\colon\lim_{n\to\omega}\|x_n\|_{2,\tau_n}=0\right\},
\end{equation*}
where $\prod_{n=1}^\infty\mathcal M_n$ consists of all sequences $(x_n)$ with $x_n\in \mathcal M_n$ and $\sup_n\|x_n\|<\infty$.
We write $(x_n)_\omega$ for the class of a bounded sequence.
The ultraproduct $\prod_\omega\mathcal M_n$ is a finite factor with faithful normal tracial state
\begin{equation*}
    \tau_\omega((x_n)_\omega)=\lim_{n\to\omega}\tau_n(x_n).
\end{equation*}
It is of type $\mathrm{II}_1$ if and only if $\dim_{\mathbb C}\mathcal M_n\to_\omega\infty$.
For a constant sequence $\mathcal M_n=\mathcal M$, we write $\mathcal M^\omega$ for the \emph{ultrapower} of $\mathcal M$.

Following Dixmier's work \cite{Dixmier-69} on central sequences, McDuff \cite[Section~4]{McDuff-70} used ultrapowers to study central sequence algebras and their connection with tensorial absorption of $\mathcal R$.
In his work on injective factors, Connes \cite[Section~V]{Connes-76} observed that $\mathcal L(\mathbb F_2)$ embeds into $\mathcal R^\omega$ and raised the question whether every separable type $\mathrm{II}_1$ factor admits such an embedding.
Ge and Hadwin \cite{GH-01} developed a systematic study of ultraproducts of operator algebras, including ultrapowers of finite von Neumann algebras.
Fang, Ge, and Li \cite[Theorem~4.5]{FGL-06} proved that $\mathcal M^\omega$ is prime for every type $\mathrm{II}_1$ factor $\mathcal M$ with separable predual.
See \cite[Sections~A.5 and~A.6]{SS-08-book} for further background on ultrapowers of type $\mathrm{II}_1$ factors and their relative commutants.
Kunnawalkam Elayavalli \cite[Problem~39]{KE-26} asked whether the ultrapower of a type $\mathrm{II}_1$ factor can ever be isomorphic to a group von Neumann algebra.
Our main theorem gives a negative answer even for embeddings and extends the conclusion to subfactors of finite Jones index.

\begin{theorem}\label{thm embedding}
Let $(\mathcal M_n,\tau_n)_{n\geqslant1}$ be a sequence of finite factors and let $\omega$ be a free ultrafilter on $\mathbb N$ such that $\dim_{\mathbb C}\mathcal M_n\to_\omega\infty$.
If $\mathcal N$ is a subfactor of finite Jones index in $\prod_\omega\mathcal M_n$, then $\mathcal N$ does not embed into $\mathcal L(G)$ for any discrete group $G$.
\end{theorem}

Taking $\mathcal N=\prod_\omega\mathcal M_n$ in \Cref{thm embedding} gives the conclusion for the ultraproduct itself.
This applies to the ultrapower of every type $\mathrm{II}_1$ factor, including $\mathcal R^\omega$, and to the matrix ultraproduct $\prod_\omega M_n(\mathbb C)$.
The main step is \Cref{prop orthogonal}, where we show that any subfactor of an ultraproduct in \Cref{thm embedding} that embeds into a group von Neumann algebra is orthogonal to a copy of $\mathcal L(\mathbb F_{\mathfrak c})$ in the ultraproduct, where $\mathfrak c$ denotes the cardinality of the continuum.

To the best of our knowledge, these ultraproducts are the first examples of type $\mathrm{II}_1$ factors shown to admit no embedding into the group von Neumann algebra of any discrete group.
The ultraproducts in \Cref{thm embedding} have nonseparable predual and the question in \cite[Remark~1.1]{CDI-23} for separable type $\mathrm{II}_1$ factors remains open.

\section{Proof of the main theorem}

Let $(\mathcal L(G),\tau)$ denote the group von Neumann algebra of a discrete group $G$.
For each $g\in G$, we denote by $u_g$ the canonical unitary operator in $\mathcal L(G)$.
Every operator $x\in\mathcal L(G)$ has a unique $L^2$-Fourier expansion
\begin{equation*}
	x=\sum_{g\in G}\widehat x(g)u_g,\qquad
	\widehat x(g)=\tau(xu_g^*),\qquad
	\sum_{g\in G}|\widehat x(g)|^2=\|x\|_2^2.
\end{equation*}
The \emph{Fourier height} \cite[(3.1)]{IPV-13} of $x\in\mathcal L(G)$ is defined as
\begin{equation*}
	h_G(x)=\sup_{g\in G}|\widehat x(g)|.
\end{equation*}
A \emph{symmetry} is a self-adjoint unitary operator.
The following technical lemma provides the symmetries used in \Cref{lem shift}.

\begin{lemma}\label{lem symmetry}
Let $\mathcal A\subseteq\mathcal L(G)$ be a diffuse abelian von Neumann subalgebra.
If $\mathcal D\subseteq\mathcal A$ is a finite-dimensional von Neumann subalgebra, then for every $\varepsilon>0$, there exists a symmetry $s\in\mathcal A$ such that $E_{\mathcal D}(s)=0$ and $h_G(s)<\varepsilon$.
\end{lemma}

\begin{proof}
Fix an integer $m>3\varepsilon^{-4}$.
Since $\mathcal D$ is finite-dimensional and abelian, $1$ is the sum of all (finitely many) minimal projections in $\mathcal D$.
Since $\mathcal A$ is diffuse, every minimal projection $p$ in $\mathcal D$ is the sum of $2m$ mutually orthogonal projections in $\mathcal A$, each of trace $\tau(p)/(2m)$.
Group these projections into pairs within each $p$, and enumerate all the pairs as $\{(e_j,f_j)\}_{j=1}^N$.
Let
\begin{equation*}
	d_j=e_j-f_j,\qquad \tau_j=\tau(e_j+f_j).
\end{equation*}
For every pair $(e_j,f_j)$ within $p$, we have $E_{\mathcal D}(e_j)=E_{\mathcal D}(f_j)=p/(2m)$.
It follows that $E_{\mathcal D}(d_j)=0$.
Moreover, $\sum_{j=1}^{N}\tau_j=1$ and $0\leqslant\tau_j\leqslant 1/m$ for $1\leqslant j\leqslant N$.

Let $\{\varepsilon_j\}_{j=1}^N$ be independent random variables taking the values $1$ and $-1$ with equal probability, and let
\begin{equation*}
	s=\sum_{j=1}^N\varepsilon_jd_j.
\end{equation*}
Then $s$ is a symmetry in $\mathcal A$ and $E_{\mathcal D}(s)=0$.
Let $\mu_g=\sum_{j=1}^N|\widehat{d_j}(g)|^2$ for each $g\in G$.
Since $|\widehat{d_j}(g)|\leqslant\|d_j\|_1=\tau_j\leqslant1/m$, we have
\begin{equation*}
	\mu_g\leqslant\sum_{j=1}^N\tau_j^2
	\leqslant\frac{1}{m}\sum_{j=1}^N\tau_j=\frac{1}{m}.
\end{equation*}
Moreover, Parseval's identity gives
\begin{equation*}
	\sum_{g\in G}\mu_g=\sum_{j=1}^N\|d_j\|_2^2
	=\sum_{j=1}^N\tau_j=1.
\end{equation*}
By the Khintchine inequality with $p=4$, we have
\begin{equation*}
	\mathbb E|\widehat s(g)|^4
	=\mathbb E\left|\sum_{j=1}^N\varepsilon_j\widehat{d_j}(g)\right|^4
	\leqslant 3\left(\sum_{j=1}^N|\widehat{d_j}(g)|^2\right)^2
	=3\mu_g^2
	\qquad(g\in G).
\end{equation*}
It follows that
\begin{equation*}
	\mathbb E\,h_G(s)^4
	\leqslant\sum_{g\in G}\mathbb E|\widehat s(g)|^4
	\leqslant 3\sum_{g\in G}\mu_g^2
	\leqslant\frac{3}{m}<\varepsilon^4.
\end{equation*}
Thus, some choice of the signs $\{\varepsilon_j\}_{j=1}^N$ gives $h_G(s)<\varepsilon$.
\end{proof}

A family $(\mathcal A_j)_{j\in J}$ of von Neumann subalgebras of an abelian tracial von Neumann algebra $(\mathcal A,\tau)$ is \emph{independent} if
\begin{equation*}
    \tau(x_1x_2\cdots x_m)=\tau(x_1)\tau(x_2)\cdots\tau(x_m)
\end{equation*}
for all distinct indices $j_1,j_2,\ldots,j_m\in J$ and all operators $x_k\in\mathcal A_{j_k}$.
We construct a sequence of independent symmetries by applying \Cref{lem symmetry} successively.

\begin{lemma}\label{lem shift}
Let $\mathcal A\subseteq\mathcal L(G)$ be a diffuse abelian von Neumann subalgebra.
Then there exists a sequence of symmetries $\{s_n\}_{n=1}^\infty$ in $\mathcal A$ such that
\begin{enumerate}
\item[$(1)$] $\tau(s_n)=0$ for $n\geqslant 1$;
\item[$(2)$] the von Neumann subalgebras $W^*(s_n)$ are independent;
\item[$(3)$] for every vector $\xi\in L^2(\mathcal L(G),\tau)$, we have
\begin{equation*}
    \lim_{\substack{\min\{m,n\}\to\infty\\m\ne n}}
    \|\xi s_n-s_m\xi\|_2=\sqrt{2}\,\|\xi\|_2.
\end{equation*}
\end{enumerate}
\end{lemma}

\begin{proof}
Let $p_1$ be a trace-half projection in $\mathcal A$.
Let $s_1=2p_1-1$ and $F_0=\varnothing$.
Suppose that $\{s_j\}_{j=1}^n$ and $\{F_j\}_{j=0}^{n-1}$ have been chosen.
Choose a nonempty finite set $F_n\subseteq G$ containing $F_{n-1}$ such that
\begin{equation*}
	x_n=\sum_{g\in F_n}\widehat{s_n}(g)u_g,
	\qquad
	\|s_n-x_n\|_2<2^{-n}.
\end{equation*}
Since the commuting symmetries $\{s_j\}_{j=1}^n$ generate a finite-dimensional von Neumann subalgebra $\mathcal D_n\subseteq\mathcal A$, by \Cref{lem symmetry}, there exists a symmetry $s_{n+1}\in\mathcal A$ such that
\begin{equation*}
	E_{\mathcal D_n}(s_{n+1})=0,\qquad
	h_G(s_{n+1})<\frac{2^{-n}}{\sqrt{|F_n|}}.
\end{equation*}
Taking traces gives $\tau(s_{n+1})=0$, which proves $(1)$.

Since $W^*(s_{n+1})=\operatorname{span}\{1,s_{n+1}\}$, we have $E_{\mathcal D_n}(x)=\tau(x)1$ for every operator $x\in W^*(s_{n+1})$.
It follows that
\begin{equation*}
	\tau(yx)=\tau\bigl(yE_{\mathcal D_n}(x)\bigr)=\tau(y)\tau(x)
	\qquad\bigl(y\in\mathcal D_n,\ x\in W^*(s_{n+1})\bigr).
\end{equation*}
Thus, $W^*(s_{n+1})$ is independent of $\mathcal D_n$, which proves $(2)$ by induction.

Let $m>n$ and $g,h\in G$.
Since $x_n$ has Fourier support in $F_n$ and $\|x_n\|_2\leqslant 1$, the Cauchy--Schwarz inequality gives
\begin{equation*}
	|\langle u_h,s_mu_gx_n\rangle|
	=|\langle u_g^*s_mu_h,x_n\rangle|
	\leqslant\sqrt{|F_n|}\,h_G(u_g^*s_mu_h)\|x_n\|_2
	\leqslant\sqrt{|F_n|}\,h_G(s_m).
\end{equation*}
Since $F_n\subseteq F_{m-1}$, the choice of $s_m$ gives $\sqrt{|F_n|}\,h_G(s_m)\leqslant 2^{-(m-1)}\leqslant 2^{-n}$.
It follows that
\begin{equation*}
	|\langle u_h,s_m u_g s_n\rangle|
	\leqslant \|s_n-x_n\|_2+\sqrt{|F_n|}\,h_G(s_m)
	\leqslant 2^{1-n}.
\end{equation*}
Similarly, $|\langle u_h,s_m u_g s_n\rangle|\leqslant 2^{1-m}$ for $m<n$.
Thus, the operator $T_{m,n}\xi=s_m\xi s_n$ is unitary on $L^2(\mathcal L(G),\tau)$ such that
\begin{equation*}
	|\langle u_h,T_{m,n}u_g\rangle|\leqslant 2^{1-\min\{m,n\}}\qquad (m\ne n).
\end{equation*}
Consequently, $T_{m,n}\to0$ in the weak-operator topology as $\min\{m,n\}\to\infty$ with $m\ne n$.
Thus, for every $\xi\in L^2(\mathcal L(G),\tau)$, we have
\begin{equation*}
	\|\xi s_n-s_m\xi\|_2^2
	=\|\xi-T_{m,n}\xi\|_2^2
	=2\|\xi\|_2^2-2\operatorname{Re}\langle\xi,T_{m,n}\xi\rangle
	\to 2\|\xi\|_2^2.
\end{equation*}
This proves $(3)$.
\end{proof}

Recall that two von Neumann subalgebras $\mathcal N$ and $\mathcal Q$ of a tracial von Neumann algebra $(\mathcal P,\tau)$ are \emph{orthogonal} if $\tau(xy)=\tau(x)\tau(y)$ for all $x\in\mathcal N$ and $y\in\mathcal Q$.
Let $\mathbb F_{\mathfrak c}$ be the free group on $\mathfrak c$ generators, where $\mathfrak c=2^{\aleph_0}$.

\begin{proposition}\label{prop orthogonal}
Suppose that $\mathcal P$ is a type $\mathrm{II}_1$ factor such that every normal trace-preserving automorphism of every separable abelian von Neumann subalgebra of $\mathcal P$ is implemented by a unitary operator in $\mathcal P$.
If a subfactor $\mathcal N\subseteq\mathcal P$ admits an embedding into $\mathcal L(G)$ for some discrete group $G$, then $\mathcal P$ contains a subfactor isomorphic to $\mathcal L(\mathbb F_{\mathfrak c})$ and orthogonal to $\mathcal N$.
\end{proposition}

\begin{proof}
Let $E_{\mathcal N}\colon\mathcal P\to\mathcal N$ be the trace-preserving conditional expectation.
Suppose first that $\mathcal N$ is of type $\mathrm{II}_1$.
Identify $\mathcal N$ with its image in $\mathcal L(G)$, where the group trace restricts to its normalized trace.
Apply \Cref{lem shift} to some diffuse abelian von Neumann subalgebra of $\mathcal N$, and let $\mathcal A=W^*(\{s_n\}_{n=1}^\infty)$ for the resulting symmetries.
Since these commuting symmetries are independent and each has trace zero, every permutation of them extends to a normal trace-preserving automorphism of $\mathcal A$.
By \cite{deBruijn-57}, there is a faithful action $w\mapsto\sigma_w$ of $\mathbb F_{\mathfrak c}$ on $\mathbb N$.
For each free generator $a\in\mathbb F_{\mathfrak c}$, the hypothesis on $\mathcal P$ gives a unitary operator $u_a\in\mathcal P$ implementing $s_n\mapsto s_{\sigma_a(n)}$ on the separable algebra $\mathcal A$.
By freeness, we obtain a unitary representation $w\mapsto u_w$ satisfying $u_w s_n=s_{\sigma_w(n)}u_w$ for all $w\in\mathbb F_{\mathfrak c}$ and $n\geqslant 1$.
To prove orthogonality, let $x_w=E_{\mathcal N}(u_w)$ for $w\in\mathbb F_{\mathfrak c}$.
Bimodularity gives
\begin{equation*}
    x_ws_n=E_{\mathcal N}(u_ws_n)
    =E_{\mathcal N}(s_{\sigma_w(n)}u_w)=s_{\sigma_w(n)}x_w
    \qquad(n\geqslant1).
\end{equation*}
Let $w\ne e$.
Since the action is faithful and $\mathbb F_{\mathfrak c}$ is torsion-free, $\sigma_w$ has infinite order and hence moves infinitely many integers.
Choose an increasing sequence $\{n_j\}_{j=1}^\infty$ of integers moved by $\sigma_w$.
Then $\sigma_w(n_j)\ne n_j$ and $\sigma_w(n_j)\to\infty$ by injectivity.
Since $x_w\in\mathcal N\subseteq\mathcal L(G)$, \Cref{lem shift}(3) gives
\begin{equation*}
    0=\lim_{j\to\infty}\|x_ws_{n_j}-s_{\sigma_w(n_j)}x_w\|_2
    =\sqrt{2}\,\|x_w\|_2.
\end{equation*}
Thus, $E_{\mathcal N}(u_w)=0$, and hence $\tau(u_w)=0$ for every $w\ne e$.
Let $\mathcal Q=W^*(\{u_w\}_{w\in\mathbb F_{\mathfrak c}})$.
Since the unitary operator $\delta_w\mapsto u_w$ from $\ell^2(\mathbb F_{\mathfrak c})$ onto $L^2(\mathcal Q)$ intertwines the left regular representation with left multiplication, it follows that $\mathcal Q\cong\mathcal L(\mathbb F_{\mathfrak c})$.
Moreover, the identity $E_{\mathcal N}(y)=\tau(y)1$ extends by normality from the span of the $u_w$ to $\mathcal Q$, proving orthogonality.

If $\mathcal N$ is a matrix factor, then $\mathcal N'\cap\mathcal P$ is of type $\mathrm{II}_1$ and contains a copy $\mathcal R_0$ of the hyperfinite type $\mathrm{II}_1$ factor.
Then $\mathcal N\vee\mathcal R_0\cong\mathcal N\mathbin{\overline\otimes}\mathcal R_0\cong\mathcal R$.
Since $\mathcal R$ is a group factor, the preceding case gives a copy of $\mathcal L(\mathbb F_{\mathfrak c})$ orthogonal to $\mathcal N\vee\mathcal R_0$, hence to $\mathcal N$.
\end{proof}

For the ultraproducts in \Cref{thm embedding}, the implementation hypothesis of \Cref{prop orthogonal} follows from unitary conjugacy of trace-preserving embeddings of separable amenable algebras \cite[Section~2.3, $1^\circ$]{Popa-14}.
For ultrapowers, see \cite[Corollary~3.4]{Sherman-09}.
The same implementation property follows from \emph{countable quantifier-free saturation}, which means that every countable family of quantifier-free conditions with parameters in $\mathcal P$ and variables in prescribed norm balls has a simultaneous solution whenever every finite subfamily can be satisfied to arbitrary accuracy.
In the proof below, the conditions simply require certain $*$-polynomials to have $2$-norm zero.
If $\mathcal N\subseteq\mathcal P$ is a subfactor of finite index $[\mathcal P:\mathcal N]=\dim_{\mathcal N}L^2(\mathcal P)<\infty$, then $L^2(\mathcal P)$ is finitely generated as a right $\mathcal N$-module.
The following corollary requires only that this module be countably generated.

\begin{corollary}\label{cor saturated}
Let $\mathcal P$ be a countably quantifier-free saturated type $\mathrm{II}_1$ factor and let $\mathcal N\subseteq\mathcal P$ be a subfactor such that $L^2(\mathcal P)$ is countably generated as a right $\mathcal N$-module.
Then $\mathcal N$ admits no embedding into a group von Neumann algebra.
In particular, the same conclusion holds if $[\mathcal P:\mathcal N]<\infty$.
\end{corollary}

\begin{proof}
Let $\alpha$ be a normal trace-preserving automorphism of a separable abelian von Neumann subalgebra $\mathcal A\subseteq\mathcal P$.
Since $\mathcal A$ is separable, we can choose projections $\{p_n\}_{n=1}^\infty$ that generate $\mathcal A$.
For a variable $u$ in the closed unit ball of $\mathcal P$, consider the following countable family of quantifier-free conditions
\begin{equation*}
    \|u^*u-1\|_2=\|uu^*-1\|_2=0,\qquad
    \|up_n-\alpha(p_n)u\|_2=0\quad(n\geqslant1).
\end{equation*}
For each $m\geqslant1$, the commuting projections $\{p_j\}_{j=1}^m$ generate a finite-dimensional abelian von Neumann algebra $\mathcal D_m$.
Since $\alpha$ is trace-preserving, there exists a unitary operator $u_m\in\mathcal P$ such that $\alpha(x)=u_mxu_m^*$ for every $x\in\mathcal D_m$.
Thus, every finite subfamily of the conditions has a solution in $\mathcal P$.
By countable quantifier-free saturation, there exists a unitary $u\in\mathcal P$ such that $up_nu^*=\alpha(p_n)$ for every $n\geqslant1$.
Since the $p_n$ generate $\mathcal A$ and both maps are normal, we obtain $uxu^*=\alpha(x)$ for every $x\in\mathcal A$.

Suppose that $\mathcal N$ admits an embedding into a group von Neumann algebra.
By \Cref{prop orthogonal}, there is a copy of $\mathcal L(\mathbb F_{\mathfrak c})$ in $\mathcal P$ orthogonal to $\mathcal N$.
Writing $u_g$ for the canonical unitary operator corresponding to $g\in \mathbb F_{\mathfrak c}$ in this copy, we have
\begin{equation*}
    E_{\mathcal N}(u_g^*u_h)=\delta_{g,h}1
    \qquad(g,h\in\mathbb F_{\mathfrak c}).
\end{equation*}
Thus, $L^2(\mathcal P)$ contains the right $\mathcal N$-module $\bigoplus_{g\in\mathbb F_{\mathfrak c}} u_gL^2(\mathcal N)$.
Since $\mathbb F_{\mathfrak c}$ is uncountable, $L^2(\mathcal P)$ is not countably generated as a right $\mathcal N$-module.
That is a contradiction.
Therefore, $\mathcal N$ admits no embedding into a group von Neumann algebra.
Moreover, $L^2(\mathcal P)$ is finitely generated as a right $\mathcal N$-module if $[\mathcal P:\mathcal N]=\dim_{\mathcal N}L^2(\mathcal P)<\infty$ by \cite{Jones-83,JS-97}.
This completes the proof.
\end{proof}

We now prove our main theorem.

\begin{proof}[Proof of \Cref{thm embedding}]
The dimension hypothesis makes $\prod_\omega\mathcal M_n$ a type $\mathrm{II}_1$ factor. It is countably saturated by
\cite[Proposition~4.11]{FHS-14-II}, and hence countably quantifier-free saturated.
The assertion follows from Corollary~\ref{cor saturated}.
\end{proof}



\begin{thebibliography}{99}






\bibitem{Boutonnet-13}
R.~Boutonnet.
\emph{$W^*$-superrigidity of mixing Gaussian actions of rigid groups.}
Adv. Math. 244 (2013), 69--90.
MR3077866


\bibitem{deBruijn-57}
N.~de Bruijn.
\emph{Embedding theorems for infinite groups.}
Nederl. Akad. Wetensch. Proc. Ser. A 60.
Indag. Math. 19 (1957), 560--569.
MR0098127


\bibitem{CDI-23}
I.~Chifan, D.~Drimbe, A.~Ioana.
\emph{Embedding universality for $\mathrm{II}_1$ factors with property~$(T)$.}
Adv. Math. 417 (2023), Paper No. 108934, 24 pp.
MR4551881


\bibitem{Connes-75}
A.~Connes.
\emph{Sur la classification des facteurs de type $\mathrm{II}$.}
C. R. Acad. Sci. Paris S\'er. A-B 281 (1975), no. 1, A13--A15.
MR0377534


\bibitem{Connes-76}
A.~Connes.
\emph{Classification of injective factors. Cases $\mathrm{II}_1$, $\mathrm{II}_\infty$, $\mathrm{III}_\lambda$, $\lambda\neq1$.}
Ann. of Math. (2) 104 (1976), no. 1, 73--115.
MR0454659


\bibitem{Dixmier-69}
J.~Dixmier.
\emph{Quelques propri\'et\'es des suites centrales dans les facteurs de type $\mathrm{II}_1$.}
Invent. Math. 7 (1969), 215--225.
MR0248534


\bibitem{FGL-06}
J.~Fang, L.~Ge, W.~Li.
\emph{Central sequence algebras of von Neumann algebras.}
Taiwanese J. Math. 10 (2006), no. 1, 187--200.
MR2186173


\bibitem{FHS-14-II}
I.~Farah, B.~Hart, D.~Sherman.
\emph{Model theory of operator algebras II: model theory.}
Israel J. Math. 201 (2014), no. 1, 477--505.
MR3265292






\bibitem{GH-01}
L.~Ge, D.~Hadwin.
\emph{Ultraproducts of $C^*$-algebras.}
Recent advances in operator theory and related topics (Szeged, 1999), 305--326.
Oper. Theory Adv. Appl., 127.
Birkh\"auser Verlag, Basel, 2001.
MR1902808








\bibitem{Ioana-11}
A.~Ioana.
\emph{$W^*$-superrigidity for Bernoulli actions of property $(T)$ groups.}
J. Amer. Math. Soc. 24 (2011), no. 4, 1175--1226.
MR2813341


\bibitem{IPV-13}
A.~Ioana, S.~Popa, S.~Vaes.
\emph{A class of superrigid group von Neumann algebras.}
Ann. of Math. (2) 178 (2013), no.~1, 231--286.
MR3043581


\bibitem{Jones-80}
V.~Jones.
\emph{A $\mathrm{II}_1$ factor anti-isomorphic to itself but without involutory antiautomorphisms.}
Math. Scand. 46 (1980), no. 1, 103--117.
MR0585235


\bibitem{Jones-83}
V.~Jones.
\emph{Index for subfactors.}
Invent. Math. 72 (1983), no. 1, 1--25.
MR0696688


\bibitem{JS-97}
V.~Jones, V.~Sunder.
\emph{Introduction to subfactors.}
London Math. Soc. Lecture Note Ser., 234.
Cambridge University Press, Cambridge, 1997. xii+162 pp.
MR1473221


\bibitem{KE-26}
S.~Kunnawalkam Elayavalli.
\emph{50 Open Problems: Ultraproduct $\mathrm{II}_1$ factors.}
(2026) arXiv:2511.20377v3.




\bibitem{McDuff-70}
D.~McDuff.
\emph{Central sequences and the hyperfinite factor.}
Proc. London Math. Soc. (3) 21 (1970), 443--461.
MR0281018


\bibitem{MvN-43}
F.~Murray, J.~von Neumann.
\emph{On rings of operators. IV.}
Ann. of Math. (2) 44 (1943), 716--808.
MR0009096




\bibitem{Ozawa-04}
N.~Ozawa.
\emph{There is no separable universal $\mathrm{II}_1$-factor.}
Proc. Amer. Math. Soc. 132 (2004), no.~2, 487--490.
MR2022373




\bibitem{Popa-14}
S.~Popa.
\emph{Independence properties in subalgebras of ultraproduct $\mathrm{II}_1$ factors.}
J. Funct. Anal. 266 (2014), no.~9, 5818--5846.
MR3182961


\bibitem{Sakai-71}
S.~Sakai.
\emph{$C^*$-Algebras and $W^*$-Algebras},
Ergeb. Math. Grenzgeb., Band 60 [Results in Mathematics and Related Areas]
Springer-Verlag, New York-Heidelberg, 1971. xii+253 pp.
MR0442701




\bibitem{Sherman-09}
D.~Sherman.
\emph{Notes on automorphisms of ultrapowers of $\mathrm{II}_1$ factors.}
Studia Math. 195 (2009), no.~3, 201--217.
MR2559172


\bibitem{SS-08-book}
A.~Sinclair, R.~Smith.
\emph{Finite von Neumann algebras and masas.}
London Math. Soc. Lecture Note Ser., 351.
Cambridge University Press, Cambridge, 2008. x+400 pp.
MR2433341



\end{thebibliography}
\end{document}